\documentclass[aap,citesort,MSNbibl,dvips]{arximspdf}

\doi{10.1214/09-AAP637}
\volume{20}
\issue{2}
\pubyear{2010}
\firstpage{367}
\lastpage{381}

\makeatletter
\renewcommand{\r}{\rightarrow}
\def\one #1{1_{\{#1\}}}
\def\eq #1{(\ref{eq:#1})}
\renewcommand{\d}{\mathrm{d}}
\newcommand{\ds}{\displaystyle}

\newtheorem{Th}{Theorem}
\makeatother

\begin{document}
\begin{frontmatter}

\title{A new formula for some linear stochastic equations with applications}
\runtitle{Linear stochastic equations}

\begin{aug}
\author[A]{\fnms{Offer} \snm{Kella}\thanksref{t1}\ead[label=e1]{Offer.Kella@huji.ac.il}\corref{}}
\and
\author[B]{\fnms{Marc} \snm{Yor}\ead[label=e2]{deaproba@proba.jussieu.fr}}
\runauthor{O. Kella and M. Yor}
\affiliation{Hebrew University of Jerusalem and Universit\'{e} Pierre et Marie Curie}
\address[A]{
Department of Statistics\\
Hebrew University of Jerusalem\\
Jerusalem 91905\\
Israel\\
\printead{e1}} 
\address[B]{
Laboratoire de Probabilit\'{e}s\\
\quad et Mod\`{e}les al\'{e}atoires\\
Universit\'{e} Pierre et Marie Curie\\
Bo\^{i}te courrier 188\\
75252 Paris Cedex 05\\
France\\
\printead{e2}}
\end{aug}
\thankstext{t1}{Supported by Grant 434/09
from the Israel Science Foundation and
the Vigevani Chair in Statistics.}
\received{\smonth{7} \syear{2009}}

\begin{abstract}
We give a representation of the solution for a stochastic linear
equation of the form $X_t=Y_t+\int_{(0,t]}X_{s-}\,\d{Z}_s$ where $Z$ is
a c\`{a}dl\`{a}g semimartingale and~$Y$ is a c\`{a}dl\`{a}g adapted
process with bounded variation on finite intervals. As an application
we study the case where $Y$ and $-Z$ are nondecreasing, jointly have
stationary increments and the jumps of $-Z$ are bounded by $1$. Special
cases of this process are shot-noise processes, growth collapse
(additive increase, multiplicative decrease) processes and clearing
processes. When $Y$ and $Z$ are, in addition, independent L\'{e}vy
processes, the resulting $X$ is called a generalized
Ornstein--Uhlenbeck process.
\end{abstract}


\end{frontmatter}

\section{Introduction}

In this paper we show that when $Z$ is a c\`
{a}dl\`{a}g adapted semimartingale and $Y$ is c\`{a}dl\`{a}g adapted
and with bounded variation on compact intervals, then the unique c\`
{a}dl\`{a}g adapted solution of $X_t=Y_t+\int_{(0,t]}X_{s-}\,\d{Z}_s$
is given via the representation $X_t=\int_{[0,t]}U_{u,t}\,\d{Y}_u$
where $U_{u,t}$ is defined by formula~\eq{Uut} below. This form seems
to be new and we note that the integral with respect to $Y$ is defined
path-wise while the integral in the integral equation can be a
stochastic integral. Of course when $Y$ is a semimartingale, one cannot
expect such a representation of the solution since $\{U_{u,t}| 0\le
u\le t\}$ is not adapted as a process indexed by $u$.

We discuss an application to the case where $Y$ and $-Z$ are
nondecreasing processes jointly having stationary increments and
subsequently specialize to cases where one or both also have
independent increments (L\'{e}vy processes). This model is a
generalization of both the shot-noise process as well as a
growth--collapse process (e.g., see, \cite{LL2008,Kella2009,GRZ2004}
and references therein) or more generally an additive increase and
multiplicative decrease process. The later have been used as models for
the TCP window size in communication networks.

We note that Jacod (\cite{Jacod79}, Theorem 6.8, page~194) and Yoeurp
and Yor \cite{YoeurpYor77} give a complete solution for the case where
the integrator is a semimartingale and the driving process is c\`{a}dl\`
{a}g, Jaschke \cite{Jaschke03} gives a derivation for the case where
the integrator does not have jumps of size $-1$, and Protter (\cite{Protter04}, Theorems 52 and 53, pages~322--323) treats the case with a
continuous integrator.

The literature related to generalized Ornstein--Uhlenbeck processes and
their applications which are directly related to some of the special
cases of the applications that we consider is huge and growing
exponentially fast. We refer the reader to \cite{BLM2006,BY2005,BBY2002,CPY1997,CPY2001,EM2005,KLM2004,Lachal2003,LM2005,LS2009,NP1996,Yor2001}
and further references therein.

\section{Main result}
With respect to some standard (right continuous augmented) filtration,
let $Y=\{Y_t| t\ge0\}$ and $Z=\{Z_t| t\ge0\}$ be two adapted c\`
{a}dl\`{a}g processes. Denote $Z_{0-}=0$, and for $t>0$, $Z_{t-}=\lim
_{s\uparrow t}Z_s$. Set $\Delta Z_t=Z_t-Z_{t-}$ when $Z$ is of bounded
variation on compact intervals (BV); set $Z^c_t=Z_t-\sum_{0\le s\le
t}\Delta Z_s$ and similarly for any other c\`{a}dl\`{a}g process
considered in this paper.
\begin{Th}\label{Th:sol}
Assume $Y$ and $Z$ are c\`{a}dl\`{a}g and
adapted, $Y$ is BV and $Z$ is a semimartingale. Then
the unique c\`{a}dl\`{a}g adapted solution to the equation
$X_t=Y_t+\int_{(0,t]}X_{s-}\,\d{Z}_s$ is
%
\begin{equation}
X_t=\int_{[0,t]}U_{u,t}\,\d{Y}_u ,
\end{equation}
where
%
\begin{eqnarray}\label{eq:Uut}
U_{u,t}=
\cases{
\ds e^{Z_t-Z_u-(1/2)([Z,Z]^c_t-[Z,Z]^c_u)}\cr
\qquad\times\ds\prod_{u<s\le t}(1+\Delta Z_s)e^{-\Delta Z_s},
 &\quad $0\le u<t$, \cr
1,&\quad $0\le u=t$
}
\end{eqnarray}
and $[Z,Z]$ is the quadratic variation process associated with $Z$.
When $Z$ is BV then \eq{Uut} reduces to
%
\begin{equation}\label{eq:UutBV}
U_{u,t}= \cases{
\ds e^{Z^c_t-Z^c_u}\prod_{u<s\le t}(1+\Delta Z_s), &\quad $0\le u<t$, \cr
1,&\quad $0\le u=t$,
}
\end{equation}
where $Z^c$ is the continuous part of $Z$ as defined earlier (rather
than the continuous martingale part of $Z$ as is customary in
stochastic calculus).
\end{Th}

\begin{pf}
Note that with $T_0=0$ and for $n\ge1$,
$T_{n}=\inf\{t>T_{n-1}| \Delta Z_t=-1\}$, then for $T_n<u\le t<T_{n+1}$
%
\begin{equation}
\frac{U_{T_n,t}}{U_{T_n,u-}}=U_{u,t}(1+\Delta Z_u) ,\qquad \frac
{U_{T_n,t}}{U_{T_n,u}}=U_{u,t} .
\end{equation}
Also, since $Y$ is a BV process, the covariation process $[Y,Z]$ is
given via $[Y,Z]_t=\sum_{0\le s\le t}\Delta Y_s\Delta Z_s$.
If one follows the solution in equation~(6.9) in Theorem~(6.8) on page~194 of
\cite{Jacod79}, then for $T_n\le t<T_{n+1}$ we have that
%
\begin{eqnarray}
X_t
&=& U_{T_n,t} \biggl(\Delta Y_{T_n}+\int
_{(T_n,t]}U^{-1}_{T_n,u-}\,\d{Y}_u-\int_{(T_n,t]}U^{-1}_{T_n,u}\,\d[Y,Z]_u
 \biggr)\nonumber\\
&=& U_{T_n,t}\Delta Y_{T_n}+\int_{(T_n,t]}U_{u,t}(1+\Delta Z_u)\,\d
{Y}_u-\sum_{T_n< u\le t}U_{u,t}\Delta Y_u\Delta Z_u\\
&=&\int
_{[T_n,t]}U_{u,t}\,\d{Y}_u,\nonumber
\end{eqnarray}
where the second equality is justified since the first integral on the
right-hand side of the first equality is a path-wise Stieltjes
integral, and the second is a sum which is also defined path-wise. If
$Y$ was a general semimartingale, then interchanging $U_{T_n,t}$ with
the integral sign like this would not be justified as the resulting
integrand would no longer be adapted.
Clearly if $n\ge1$, then $U_{u,t}=0$ for $u<T_n$, and thus
%
\begin{equation}
X_t=\int_{[0,t]}U_{u,t}\,\d{Y}_u .
\end{equation}
Since this holds for all $n$, the proof for the more general case is
complete. For the case where $Z$ is BV, it is evident that $[Z,Z]^c=0$,
and it is easy to check that $\sum_{u<s\le t}{\Delta Z_s}$ is
convergent (actually, absolutely convergent), and hence the result
follows.
\end{pf}

Of course one may also define the counting process,
%
\begin{equation}
N_t=\sum_{0<s\le t}\one{\Delta Z_s=-1},
\end{equation}
which is a.s. finite for all $t\ge0$ and right-continuous (possibly
a.s. identically zero or terminating), and write

\begin{equation}
X_t=\int_{[T_{N_t},t]}U_{u,t}\,\d{Y}_u .
\end{equation}

It is worth while to note that for the case where $Z$ is also a BV
process, there is a more direct proof involving (path-wise) Stieltjes
integration which can be taught in a classroom as follows. Write
$Z=A-B$, where $A$ and $B$ are right-continuous and nondecreasing and
have no jump points in common. Write $A^d_t=A_t-A^c_t=\sum_{0<s\le
t}\max(\Delta Z_s,0)$ and similarly for $B$. Observe that by right continuity
$\Delta A_t$, $\Delta B_t$, $A^d_t-A_0$ and $B^d_t-B_0$ all converge to
zero as $t\downarrow0$. In particular, for every $t$ for which $-1\le
\Delta B_s\,(\le0)$ for $0<s\le t$, we have that
%
\begin{equation}
1+A^d_t-A_0\le\prod_{0<s\le t}(1+\Delta A_s)\le e^{A^d_t-A_0}
\end{equation}
and
%
\begin{equation}
1+B^d_t-B_0\le\prod_{0<s\le t}(1+\Delta B_s)\le e^{B^d_t-B_0}
\end{equation}
which implies that
%
\begin{equation}
\prod_{0<s\le t}(1+\Delta Z_s)= \biggl(\prod_{0<s\le t}(1+\Delta
A_s) \biggr) \biggl(\prod_{0<s\le t}(1+\Delta B_s) \biggr)\r{1}
\end{equation}
as $t\downarrow0$.

Now note that with $C_t=e^{Z^c_t}$ and $D_t=\prod_{0<s\le t}(1+\Delta
Z_s)$, ordinary (Stieltjes) integration by parts yields
%
\begin{equation}\hspace*{22pt}
U_t\equiv C_tD_t=C_{0+}D_{0+}+\int_{(0,t]}D_{s-}\,\d{C}_s+\int
_{(0,t]}C_{s-}\,\d{D}_s+\sum_{0<s\le t}\Delta C_s\Delta D_s
,
\end{equation}
and it is easy to check that the continuity of $C$ and the fact that
$\d{C}_t=C_t \,\d{Z}^c_t$ imply that
%
\begin{equation}\label{eq:U}
U_t=1+\int_{(0,t]}U_{s-}\,\d{Z}_s .
\end{equation}
With this formula established, it is clear that if we denote $U_{u,t}$
as in \eq{UutBV}, then in an identical way to which \eq{U} was
obtained we have (path-wise) that
%
\begin{equation}
U_{u,t}=1+\int_{(u,t]}U_{u,s-}\,\d{Z}_s
\end{equation}
for all $0\le u\le t$.

Now, if $X_t=\int_{[0,t]}U_{s,t}\,\d{Y}_s$, then $X_{t-}=\int
_{[0,t)}U_{s,t-}\,\d Y_s$ and thus $\int_{(0,t]}X_{s-}\,\d{Z}_s$ is given by
\begin{eqnarray}
\int_{(0,t]}\int_{[0,s)}U_{u,s-}\,\d{Y}_u \,\d{Z}_s
&=&\int_{[0,t)}\int
_{(u,t]}U_{u,s-}\,\d{Z}_s\,\d{Y}_u\nonumber\\[-8pt]\\[-8pt]
&=&\int_{[0,t)}(U_{u,t}-1)\,\d{Y}_u,\nonumber
\end{eqnarray}
but since $U_{t,t}=1$ we can include $t$ in the domain of integration
without changing the value which gives
%
\begin{equation}
\int_{(0,t]}X_{s-}\,\d{Z}_s=\int_{[0,t]}(U_{u,t}-1)\,\d{Y}_u=X_t-Y_t
\end{equation}
as required.

\section{Applications}
Assume that $Y$ and $Z$ are right-continuous and nondecreasing jointly
having stationary increments in the strong sense that the law of
$\theta_s (Y,Z)$ is independent of $s$ where
%
\begin{equation}
\theta_s(Y(t),Z(t))=\bigl(Y(t+s)-Y(s),Z(t+s)-Z(s)\bigr) .
\end{equation}
It is standard to (uniquely) extend $(Y,Z)$ to be a double sided
process having stationary increments, that is, that $t\in\mathbb{R}$
rather than $t\ge0$, thus we assume it at the outset. Finally we
assume that $Z$ has jumps bounded by $1$. Without loss of generality
let us assume that $Y_0=Z_0=0$, otherwise we perform what follows for
$Y-Y_0$ and $Z-Z_0$ which also have stationary increments. We consider
the unique process~$X$ defined via $X_t=X_0+Y_t-\int_{(0,t]}X_{s-}\,\d
{Z}_s$ for $t\ge0$ where $X_0$ is almost surely finite; the unique
solution of which is
%
\begin{equation}\label{eq:sollevy}\hspace*{10pt}
X_t=X_0e^{-Z^c_t}\prod_{0<s\le t}(1-\Delta Z_s)+\int
_{(0,t]}e^{-(Z^c_t-Z^c_u)}\prod_{u<s\le t}(1-\Delta Z_s)\,\d{Y}_u,
\end{equation}
where an empty product (when $u=t$ or when $t=0$ on the right) is
defined to be~$1$.

Special cases of such processes are the shot-noise processes in which
$Z_t=rt$ and $Y$ are compound Poisson, growth collapse or additive
increase multiplicative decrease (AIMD) processes in which $Y_t=rt$ and
usually $Z=qN_\lambda$ where $N_\lambda$ is a Poisson process with
rate $\lambda$, and $0<q<1$, as well as clearing processes where $Z$
is a Poisson process or, more generally, a renewal counting process
(see, e.g., \cite{KPS2003,Kella1998}).

Consider the nondecreasing processes
%
\begin{equation}
J_t=Z^c_t-\sum_{0<s\le t}\log(1-\Delta Z_s)\one{\Delta Z_s<1}
,
\end{equation}
and $N_t=\sum_{0<s\le t}\one{\Delta Z_s=1}$. Then it is clear that
$Y,J,N$ jointly have stationary increments (in the strong sense), and
from \eq{sollevy} we have
%
\begin{equation}
X_t=X_0e^{-J_t}\one{N_t=0}+\int_{(0,t]}e^{-(J_t-J_s)}\one
{N_t-N_s=0}\,\d{Y}_s .
\end{equation}
If $\int_{(-\infty,0]}e^{J_s}\,\d{Y}_s$ is a.s. finite (recalling that
for $s\le0$, $J_s\le J_0=0$), then setting $X^*_t=\int_{(-\infty
,t]}e^{-(J_t-J_s)}\one{N_t-N_s=0}\,\d{Y}_s$ it is clear\vspace*{1pt} that $X^*$ is a
stationary process. Moreover, if, in addition, either $\lim_{t\r
\infty}N_t\ge1$ a.s. (equivalently, $T_1=\inf\{t| \Delta Z_t=1\}$
is a.s. finite) or $J_t\r\infty$ a.s. as $t\r\infty$, then
$|X^*_t-X_t|\r{0}$ a.s. as $t\r\infty$, and thus for any a.s. finite
initial $X_0$, a limiting distribution exists which is distributed like $X^*_0$.

In fact, when $X_0$ is independent of $(Y,Z)$, then shifting by $-t$,
noting that $\theta_{-t} J_s=J_{s-t}-J_{-t}$ (so that $\theta
_{-t}J_t=0$) and similarly for $N$ and $Y$, it is clear that~$X_t$ has
the same distribution as
\begin{eqnarray}\label{eq:reverse}
&& X_0e^{J_{-t}}\one{N_{-t}=0}+\int_{(0,t]}e^{J_{s-t}}\one
{N_{s-t}=0}\,\d{Y}_{s-t}\nonumber\\[-8pt]\\[-8pt]
&&\qquad= X_0e^{J_{-t}}\one{N_{-t}=0}+\int_{(-t,0]}e^{J_s}\one{N_s=0}\,\d
{Y}_s .\nonumber
\end{eqnarray}
In particular, this implies that when $X_0=0$, then $X_t$ is
stochastically increasing in $t\ge0$.

Let us summarize our findings as follows.
\begin{Th}
If $\int_{(-\infty,0]}e^{J_s}\,\d{Y}_s<\infty$ a.s., and either
$T_1<\infty$ a.s. or $J_t\r\infty$ a.s. as $t\r\infty$, then $X$
has the unique stationary version
%
\begin{equation}
X^*_t=\int_{(-\infty,t]}e^{-(J_t-J_s)}\one{N_t-N_s=0}\,\d{Y}_s
,
\end{equation}
and for every initial a.s. finite $X_0$, $X_t$ converges in
distribution to $X^*_0$. Moreover, when $X_0=0$ a.s., then $X_t$ is
stochastically increasing in $t\ge0$.
\end{Th}

We note that when $(Y,Z)$ also have independent increments so that they
form a L\'{e}vy process, then the negative of the time reversed process
is a left-continuous version of the forward process, and thus in this
case [when $X_0$ is independent of $(Y,Z)$], $X_t$ is also distributed like
%
\begin{equation}\label{eq:ou}
X_0e^{-J_{t}}\one{N_{t}=0}+\int_{(0,t]}e^{-J_s}\one{N_s=0}\,\d{Y}_s
\end{equation}
which is also the consequence of the usual time reversal argument for
L\'{e}vy processes. In what follows we will consider special cases of
this structure.

We observe that in the general case $N$ is a simple (i.e., a.s. $\Delta
N_t\in\{0,1\}$ for all $t$) counting process associated with a time
stationary point process. Special cases of such processes are Poisson
processes and delayed renewal processes where the delay has the
stationary excess lifetime distribution associated with the subsequent
i.i.d. inter-renewal times. We will consider this special case a bit
later.\looseness=1

\subsection{$EX_t$ for independent $X_0$, $Y$, $Z$}
Since $Y$ has stationary increments, it follows that $EY_t=EY_1t$. From
\eq{reverse} we have that when $EY_1$ and $EX_0$ are finite, then for
$t\ge0$,
%
\begin{equation}
EX_t=EX_0Ee^{J_{-t}}\one{N_{-t}=0}+EY_1\int_{-t}^0Ee^{J_s}\one
{N_s=0}\,\d{s}
,
\end{equation}
and since for $s\le0$, we have that $J_s=-(J_0-J_s)$ is distributed
like $-J_{-s}=-(J_{-s}-J_0)$, and similarly for $N$, we have that
%
\begin{equation}\label{eq:meanind}
EX_t=EX_0Ee^{-J_{t}}\one{N_{t}=0}+EY_1\int_0^tEe^{-J_s}\one{N_s=0}\,\d
{s} .
\end{equation}

\subsection{\texorpdfstring{$EX_t$ for independent
$X_0$, $Y$, $Z$ with L\'{e}vy $Z$}{$EX_t$ for independent $X_0$, $Y$, $Z$ with Levy $Z$}}

Here $Z$ is a subordinator with Laplace--Stieltjes exponent
$-\eta_z(\alpha)=\log Ee^{-\alpha Z_1}$ where, for $\alpha\ge0$,
%
\begin{equation}
\eta_z(\alpha)=c_z\alpha+\int_{(0,1]}(1-e^{-\alpha x})\nu_z(\d{x})
\end{equation}
with $c_z\ge0$ and $\int_{(0,1]}x\nu_z(\d{x})<\infty$. Since the
jumps of $Z$ are bounded above by~$1$, then $\nu_z((1,\infty))=0$.

In this case $Z^c_t=c_zt$, $N$ is a Poisson process with rate $\lambda
=\nu_z\{1\}$ which is independent of the subordinator,
%
\begin{equation}
J_t=c_zt-\sum_{0<s\le t}\log(1-\Delta Z_s)\one{\Delta Z_s<1}
;
\end{equation}
the L\'{e}vy measure of which, call it $\nu_j$, is defined via $\nu
_j((a,b])=\nu_z((1-e^{-a},1-e^{-b}])$ for $0<a<b<\infty$ and with exponent
\begin{eqnarray}
\eta_j(\alpha)
&=&c_z\alpha+\int_{(0,\infty)}(1-e^{-\alpha x})\nu
_j(\d{x})\nonumber\\[-8pt]\\[-8pt]
&=&c_z\alpha+\int_{(0,1)}\bigl(1-(1-x)^\alpha\bigr)\nu_z(\d{x}) ,\nonumber
\end{eqnarray}
so that for $\alpha>0$,
%
\begin{equation}\label{eq:withlambda}
\eta_j(\alpha)+\lambda=c_z\alpha+\int_{(0,1]}\bigl(1-(1-x)^\alpha\bigr)\nu
_z(\d{x}) .
\end{equation}
We note that
%
\begin{equation}
\int_{(0,\infty)}\min(x,1)\nu_j(\d{x})=\int_{(0,1)}\min\bigl(-\log
(1-x),1\bigr)\nu_z(\d{x}), \end{equation}
and since $-\log(1-x)\le\frac{x}{1-x}\le xe$ for $0<x\le1-e^{-1}$,
the right-hand side is dominated above by $e\int_{(0,1)}x\nu_z(\d
{x})<\infty$, so that $\nu_j$ is indeed the proper L\'{e}vy measure
of a subordinator.
Now, for this case, $Ee^{-J_s}=e^{-\eta_j(1)s}$ where
\begin{eqnarray}
\eta_j(1)
&=&c_z+\int_{(0,1)}\bigl(1-(1-x)^1\bigr)\nu_z(\d{x})\nonumber\\[-8pt]\\[-8pt]
&=&c_z+\int
_{(0,1)}x\nu_z(\d{x})
=\eta_z'(0)-\lambda\nonumber
\end{eqnarray}
recalling $\lambda=\nu_z\{1\}$. Therefore, $Ee^{-J_s}\one
{N_s=0}=e^{-(\eta'_z(0)-\lambda)s}e^{-\lambda s}=e^{-\eta'_z(0)s}$
so that in this case, since $\eta'_z(0)=c_z+\int_{(0,1]}x\nu_z(\d
{x})=EZ_1$, \eq{meanind} becomes
%
\begin{equation}
EX_t=EX_0 e^{-EZ_1t}+\frac{EY_1}{EZ_1}(1-e^{-EZ_1t}) .
\end{equation}
Recall that here $Y$ need not have independent increments.

\subsection{\texorpdfstring{Independent $X_0$, $Y$, $Z$ with L\'{e}vy $Y$}{Independent $X_0$, $Y$, $Z$ with Levy $Y$}}

Since for every $0=t_0<t_1<\cdots<t_n=t$ the independence between $Y$
and $Z$ and hence the independence of $Y$ and~$J$, yield
\begin{eqnarray}
&& E \Biggl[\exp \Biggl(-\alpha\sum_{i=1}^n e^{-J_{t_{i-1}}}\one
{N_{t_{i-1}}=0}(Y_{t_i}-Y_{t_{i-1}}) \Biggr) \Big|Z
\Biggr]\nonumber\\[-8pt]\\[-8pt]
&&\qquad
=\prod_{i=1}^n\exp \bigl(
-\eta_y\bigl(\alpha e^{-J_{t_{i-1}}}\one
{N_{t_{i-1}}=0}\bigr)(t_i-t_{i-1}) \bigr).\nonumber
\end{eqnarray}
It thus follows, as in equation~(5.9) of \cite{KW1999} for the more
general multivariate case and in Proposition~1 of \cite{NP1996} for
the case where $Y$ and $Z$ are compound Poisson, that
%
\begin{eqnarray}\label{eq:inteZdY}
&&E \biggl[\exp \biggl( -\alpha\int_{(0,t]}e^{-J_s}\one{N_s=0}\,\d
{Y}_s \biggr) \Big| Z \biggr]\nonumber\\[-8pt]\\[-8pt]
&&\qquad=\exp \biggl(-\int_0^t\eta_y\bigl(\alpha
e^{-J_s}\one{N_s=0}\bigr)\,\d{s} \biggr) .\nonumber
\end{eqnarray}
This implies, as in Theorem~5.1 of \cite{KW1999}, that the conditional
distribution of $\int_{(0,t]}e^{-J_s}\one{N_s=0}\,\d{Y}_s$ given $Z$
is infinitely divisible, as on the right-hand side, $-\eta_y/n$ is
also a Laplace--Stieltjes exponent of a subordinator.

Equation~\eq{inteZdY}, with $\xi_0(\alpha)=Ee^{-\alpha X_0}$,
$a\wedge b=\min(a,b)$, and recalling
%
\begin{equation}T_1=\inf\{t| \Delta Z_t=1\}=\inf\{t| N_t>0\}
\end{equation}
yields
%
\begin{eqnarray}\label{eq:LSTX_t}
Ee^{-\alpha X_t}
&=& E\xi_0 \bigl(\alpha e^{-J_t}\one{N_t=0}
\bigr)\exp \biggl(-\int_0^{t}\eta_y(\alpha e^{-J_s})\one{N_s=0}\,\d
{s} \biggr)
\nonumber\\
&=&
E\xi_0 (\alpha e^{-J_t} )\exp \biggl(-\int
_0^{t}\eta_y(\alpha e^{-J_s})\,\d{s} \biggr)\one{T_1>t}\\
&&{} +
E\exp \biggl(-\int_0^{T_1}\eta_y(\alpha e^{-J_s})\,\d{s} \biggr)\one
{T_1\le t}
 .\nonumber
\end{eqnarray}
Clearly, when either $T_1<\infty$ a.s. or $J_t\r\infty$ a.s. as $t\r
\infty$, then
%
\begin{equation}
\lim_{t\r\infty}Ee^{-\alpha X_t}=E\exp \biggl(-\int_0^{T_1}\eta
_y(\alpha e^{-J_s})\,\d{s} \biggr).
\end{equation}
We now observe that if $N$ and $J$ are independent, as for instance in
the case where~$Z$ is a subordinator, and $N$ is the counting process
associated with a time stationary version of a renewal process the
latter having inter-renewal time distribution $F$ having a finite mean
$\mu$, then it is well known that $N$ is a delayed renewal process in
which the times between the $(i-1)$th and $i$th jumps are distributed $F$
for $i\ge2$ and the time until the first jump (i.e., the delay) has a
distribution with density $f_e(t)=(1-F(t))/\mu$. Therefore, in this case,
%
\begin{equation}\qquad
E\exp \biggl(-\int_0^{T_1}\eta_y(\alpha e^{-J_s})\,\d{s} \biggr)=\int
_0^\infty E\exp \biggl(-\int_0^{t}\eta_y(\alpha e^{-J_s})\,\d{s}
\biggr)f_e(t)\,\d{t} .
\end{equation}

Differentiating the right-hand side of the first equality in \eq
{LSTX_t} once and setting $\alpha=0$ gives \eq{meanind} as expected,
while for the case where $X_0=0$ a.s., differentiating twice and
setting $\alpha=0$ yields
%
\begin{equation}\quad \ \ \ \ 
EX_t^2=(\eta'_y(0))^2E \biggl(\int_0^te^{-J_s}\one{N_s=0}\,\d{s}
\biggr)^2-\eta''_y(0)E\int_0^te^{-2J_s}\one{N_s=0}\,\d{s} .
\end{equation}

\subsection{\texorpdfstring{$EX_t^2$ for independent $Y$, $Z$ with
L\'{e}vy $Y,Z$ and $X_0=0$}{$EX_t^2$ for independent $Y$, $Z$ with Levy $Y,Z$ and $X_0=0$}}

We note that for every $\beta>0$, $E\int_0^t e^{-\beta J_s}\one
{N_s=0}\,\d{s}=\frac{1-e^{-(\eta_j(\beta)+\lambda)t}}{\eta_j(\beta
)+\lambda}$, where $\lambda=\nu_z\{1\}$. Also, note that since
$N_u\le N_s$ for $u\le s$,
%
\begin{eqnarray}
  \biggl(\int_0^te^{-J_s}\one{N_s=0}\,\d{s} \biggr)^2
  &=& 2\int
_0^t\int_0^s e^{-J_s-J_u}\one{N_s=0}\,\d{u}\,\d{s}\nonumber\\[-8pt]\\[-8pt]
&=&
2\int_0^t \int_0^s e^{-(J_s-J_u)}e^{-2J_u}\one{N_s=0}\,\d{u}\,\d{s},\nonumber
\end{eqnarray}
and therefore (using Fubini and the stationary independent increments
property of~$J$), the expected value of the left-hand side is
%
\begin{eqnarray}
&&\hspace*{20pt}
2\int_0^t\int_0^s e^{-(\eta_j(1)+\lambda)(s-u)}e^{-(\eta
_j(2)+\lambda)u}\,\d{u}\,\d{s}\nonumber\\[-8pt]\\[-8pt]
&&\qquad\hspace*{20pt}
=2\frac{(1-e^{-(\eta_j(1)+\lambda)t})/(\eta_j(1)+\lambda)-(1-e^{-(\eta_j(2)+\lambda)t})/(\eta
_j(2)+\lambda)}{\eta_j(2)-\eta_j(1)} .\nonumber
\end{eqnarray}
Finally, we observe that for every positive integer $n$, we obtain
[recall \eq{withlambda}]
%
\begin{eqnarray}
\eta_j(n)+\lambda
&=&c_zn+\int_{(0,1]}\bigl(1-(1-x)^n\bigr)\nu_z(\d{x})\nonumber\\[-8pt]\\[-8pt]
&=&c_z
n+\sum_{k=1}^n\pmatrix{n\cr k}(-1)^{k-1}\int_{(0,1]}x^k\nu_z(\d{x}),\nonumber
\end{eqnarray}
and since, $\eta_z^{(0)}(0)=\eta_z(0)=0$, $\eta_z'(0)=c_z+\int
_{(0,1)}x\nu_z(\d{x})$ and $\eta_z^{(k)}(0)=\break(-1)^{k-1}\int
_{(0,1]}x^k\nu_z(\d{x})$, for $k\ge2$, it holds that
%
\begin{equation}\label{eq:binomial}
\eta_j(n)+\lambda=\sum_{k=0}^n\pmatrix{n\cr k}\eta_z^{(k)}(0) .
\end{equation}
In particular $\eta_j(1)+\lambda=\eta'_z(0)=c_z+\int_{(0,1]}x\nu
(\d{x})$ and $\eta_j(2)+\lambda=2\eta'_z(0)+\eta''_z(0)$, so that
$\eta_j(2)-\eta_j(1)=\eta'_z(0)+\eta''_z(0)$.

To summarize, when $EX_0=0$, we have
%
\begin{eqnarray}\label{eq:2ndmoment}
\hspace*{27pt}EX_t^2
&=&
2(\eta'_y(0))^2\nonumber\\
&&{}\times\bigl(\bigl(1-e^{-\eta'_z(0)t}\bigr)/\eta'_z(0)
-
\bigl(1-e^{-\bigl(2\eta'_z(0)+\eta''_z(0)\bigr)t}\bigr)/\bigl(2\eta'_z(0)+\eta
''_z(0)\bigr)\bigr)\nonumber\\[-8pt]\\[-8pt]
&&\hspace*{14pt}/\bigl(\eta'_z(0)+\eta''_z(0)\bigr)\nonumber\\
&&{}-\eta''_y(0)\frac{1-e^{-(2\eta'_z(0)+\eta''_z(0))t}}{2\eta
'_z(0)+\eta''_z(0)}\nonumber
\end{eqnarray}
which converges to
%
\begin{equation}\label{eq:2ndmoment}
\frac{2(\eta'_y(0))^2-\eta'_z(0)\eta''_y(0)}{\eta'_z(0)(2\eta
'_z(0)+\eta''_z(0))}=
\frac{ (\eta'_y(0)/\eta'_z(0))^2-\eta
''_y(0)/(2\eta'_z(0))}
{1+\eta''_z(0)/(2\eta'_z(0))}
\end{equation}
as $t\r\infty$. We note that as $\nu_z(1,\infty)=0$, then clearly
whenever either $c_z>0$ or $\nu_z(0,1)\not= 0$ (i.e., $Z-N$ is not
identically zero), it holds that
%
\begin{equation}\eta'_z(0)=c_z+\int_{(0,1]}x\nu_z(\d{x})>\int
_{(0,1]}x^2\nu_z(\d{x})=-\eta''_z(0) .
\end{equation}

\subsection{\texorpdfstring{L\'{e}vy $Z$, linear $Y$ and $X_0=x$}{Levy $Z$, linear $Y$ and $X_0=x$}}
It is of interest to consider the special case where $Y_t=rt$ for some
$r>0$ and $X_0=x$ for some $x\ge0$. For the case where $Z$ is compound
Poisson this model becomes the growth--collapse process from \cite{LL2008} where the computation of transient moments turns out to be
especially tractable. Since
%
\begin{equation}
\frac{X_t}{r}=\frac{x}{r}+t-\int_{(0,t]}\frac{X_{s-}}{r}\,\d{Z}_s
\end{equation}
we may without loss of generality assume that $r=1$. Recall \eq{ou}.
Following the ideas in the proof of Proposition~3.1 of \cite{CPY1997}, we
first write for $a\ge0$ and $b\ge1$,
%
\begin{eqnarray}
 Ee^{-aJ_t} \biggl(\int_0^te^{-J_s}\,\d{s} \biggr)^b
&=& bEe^{-aJ_t}\int_0^t \biggl(\int_u^te^{-J_s}\,\d{s}
\biggr)^{b-1}e^{-J_u}\,\d{u}\nonumber\\
&=& b\int_0^tEe^{-a(J_t-J_u)} \biggl(\int_u^te^{-(J_s-J_u)}\,\d
{s} \biggr)^{b-1}e^{-(a+b)J_u}\,\d{u}\\
&=& b\int_0^t e^{-\eta_j(a+b)u}Ee^{-aJ_{t-u}} \biggl(\int
_0^{t-u}e^{-J_s}\,\d{s} \biggr)^{b-1} .\nonumber
\end{eqnarray}
Thus, if $T\sim\exp(\theta)$ for some $\theta>0$ and is independent
of $Z$, then since the conditional distribution of $T-u$ given $T>u$ is
the same as that of $T$ (memoryless property), it readily follows that
%
\begin{equation}
Ee^{-aJ_T} \biggl(\int_0^Te^{-J_s}\,\d{s} \biggr)^b
=\frac{b}{\eta_j(a+b)+\theta}Ee^{-aJ_T} \biggl(\int_0^Te^{-J_s}\,\d
{s} \biggr)^{b-1}.
\end{equation}
For $a=0$ we have that, since $T_1\wedge T\sim\exp(\lambda+\theta)$
and $\int_0^Te^{-J_s}\one{N_s=0}\,\d{s}=\int_0^{T_1\wedge
T}e^{-J_s}\,\d{s}$,
%
\begin{equation}\label{eq:aeq0bge1}\quad
E \biggl(\int_0^Te^{-J_s}\one{N_s=0}\,\d{s} \biggr)^b
=\frac{b}{\eta_j(b)+\lambda+\theta}E \biggl(\int_0^Te^{-J_s}\one
{N_s=0}\,\d{s} \biggr)^{b-1} .
\end{equation}
For $a>0$ we have, from the fact that $T_1\wedge T$ is independent of
$\one{T_1>T}$, that
%
\begin{eqnarray}
&& Ee^{-aJ_T}\one{N_T=0} \biggl(\int_0^Te^{-J_s}\one{N_s=0}\,\d
{s} \biggr)^b\nonumber\\
&&\qquad= Ee^{-aJ_{T_1\wedge T}}\one{T_1>T} \biggl(\int
_0^{T_1\wedge T}e^{-J_s}\,\d{s} \biggr)^b\\
&&\qquad=  \frac{\theta
}{\lambda+\theta}Ee^{-aJ_{T_1\wedge T}} \biggl(\int_0^{T_1\wedge
T}e^{-J_s}\,\d{s} \biggr)^b\nonumber
\end{eqnarray}
and thus
%
\begin{eqnarray}\label{eq:agt0bge1}
&& Ee^{-aJ_T}\one{N_T=0} \biggl(\int_0^Te^{-J_s}\one{N_s=0}\,\d
{s} \biggr)^b\nonumber\\
&&\qquad= Ee^{-aJ_T}\one{N_T=0} \biggl(\int_0^Te^{-J_s}\,\d
{s} \biggr)^b\\
&&\qquad= \frac{b}{\eta_j(a+b)+\lambda+\theta}Ee^{-aJ_T}\one
{N_T=0} \biggl(\int_0^Te^{-J_s}\,\d{s} \biggr)^{b-1} .\nonumber
\end{eqnarray}
Clearly, when $b=0$ and $a>0$ we have that
%
\begin{equation}\label{eq:agt0beq0}
Ee^{-aJ_T}\one{N_T=0}=e^{-(\eta_j(a)+\lambda)T}=\frac{\theta}{\eta
_j(a)+\lambda+\theta} .
\end{equation}
Now
%
\begin{eqnarray}\label{eq:nmoment}
 EX_T^n
 &=& E \biggl(xe^{-J_T}\one{N_T=0}+\int_0^Te^{-J_s}\one
{N_s=0}\,\d{s} \biggr)^n\nonumber\\
 &=& \sum_{k=1}^{n}\pmatrix{n\cr k}x^kEe^{-kJ_T}\one{N_T=0}
\biggl(\int_0^Te^{-J_s}\,\d{s} \biggr)^{n-k}\\
&&{}+E \biggl(\int_0^Te^{-J_s}\one{N_s=0}\,\d{s} \biggr)^n,\nonumber
\end{eqnarray}
and denoting [recall \eq{binomial}]
%
\begin{equation}\hspace*{10pt}\mu_i=\eta_j(i)+\lambda=c_zi+\int
_{(0,1]}\bigl(1-(1-x)^i\bigr)\nu_z(\d x)=\sum_{k=0}^i\pmatrix{i\cr k}\eta
^{(k)}_z(0) ,
\end{equation}
it follows from \eq{aeq0bge1}, \eq{agt0bge1}, \eq{agt0beq0} and \eq{nmoment}, with some manipulations, that
%
\begin{eqnarray}
EX_T^n
&=& \frac{n!}{\prod_{i=1}^n\mu_i} \Biggl(\sum
_{k=1}^{n}\frac{x^k\prod_{i=1}^k\mu_i}{k!} \Biggl(\prod
_{i=k+1}^n\frac{\mu_i}{\mu_i+\theta}-
\prod_{i=k}^n\frac{\mu_i}{\mu_i+\theta} \Biggr) \nonumber\\[-8pt]\\[-8pt]
&&\hspace*{183.5pt}{}+\prod_{i=1}^n\frac{\mu_i}{\mu_i+\theta} \Biggr),\nonumber
\end{eqnarray}
where an empty product is defined to be $1$.
Finally, noting that $EX_T^n=\break\int_0^\infty e^{-\theta t}\,\d{E}X_t^n$
it follows that if $\{E_i| i\ge1\}$ are i.i.d. random variables with
distribution $\exp(1)$, then $E_i/\mu_i\sim\exp(\mu_i)$. It is
well known and easy to check that
%
\begin{equation}
\prod_{i=k}^n\frac{\mu_i}{\mu_i+\theta}=\int_0^\infty e^{-\theta
t}\,\d{P} \Biggl[\sum_{i=k}^n\frac{E_i}{\mu_i}\le t \Biggr];
\end{equation}
hence, for $1\le k\le n$,
%
\begin{equation}
\prod_{i=k+1}^n\frac{\mu_i}{\mu_i+\theta}-\prod_{i=k}^n\frac{\mu
_i}{\mu_i+\theta}=\int_0^\infty e^{-\theta t}
\,\d{P} \Biggl[\sum_{i=k+1}^n\frac{E_i}{\mu_i}\le t<\sum_{i=k}^n\frac
{E_i}{\mu_i} \Biggr] ,
\end{equation}
and thus we have the following somewhat curious result.
\begin{Th}\label{th:puredeath}
Let $p_{ij}(t)$ be the transition matrix function of a pure death
process $D=\{D_t| t\ge0\}$ with death rates $\mu_i$, $i\ge1$ ($0$
is absorbing). Then
%
\begin{eqnarray}
EX_t^n&=& \frac{n!}{\prod_{i=1}^n\mu_i} \Biggl(p_{n0}(t)+\sum
_{k=1}^{n}\frac{x^k\prod_{i=1}^k\mu_i}{k!}p_{nk}(t) \Biggr)\nonumber\\[-8pt]\\[-8pt]
 &=& \frac{n!}{\prod_{i=1}^n\mu_i}E \Biggl[\prod
_{i=1}^{D_t}\frac{x\mu_i}{i}\Big| D_0=n \Biggr],\nonumber
\end{eqnarray}
where an empty product is $1$.
\end{Th}

In particular, when $x=0$, then
%
\begin{eqnarray}\label{eq:xeq0}
EX_t^n
&=& \frac{n!}{\prod_{i=1}^n\mu_i}p_{n0}(t)\nonumber\\
&=&
n!\mathop{\mathop{\int\cdots\int}_{\sum_{i=1}^nx_i\le t}}_ {x_1,\ldots
,x_n\ge0}\exp \Biggl(-\sum_{i=1}^n\mu_ix_i \Biggr)\,\d{x}_1\cdots\,\d
{x}_n \\
&=& n!t^n\mathop{\mathop{\int\cdots\int}_{\sum
_{i=1}^nx_i\le1}}_{ x_1,\ldots,x_n\ge0}\exp \Biggl(-t\sum
_{i=1}^n\mu_ix_i \Biggr)\,\d{x}_1\cdots\,\d{x}_n.\nonumber
\end{eqnarray}

In fact, one may also give a finite simple algorithm with which to
compute $EX_t^n$. For the sake of brevity we do it only for the case
$x=0$. This can be done similarly to the Brownian motion in the proof
of Theorem~1 on page~31 of \cite{Yor2001} or, equivalently, directly
from \eq{xeq0} as follows.
Set $f_0=0$ and for $n\ge1$ and $0<a_1<a_2<\cdots<a_n$, let
%
\begin{eqnarray}
f_n(a_1,\ldots,a_n)
&=& \mathop{\mathop{\int\cdots\int}_{\sum
_{i=1}^nx_i\le1}}_{ x_1,\ldots,x_n\ge0}\exp \Biggl(-\sum
_{i=1}^na_ix_i \Biggr)\,\d{x}_1\cdots\,\d{x}_n\nonumber\\
&=&
\mathop{\mathop{\int\cdots\int}_{\sum_{i=2}^{n}x_i\le1}}_{ x_2,\ldots
,x_{n}\ge0}
 \biggl(\int_0^{1-\sum_{i=2}^{n}x_i}e^{-a_1x_1}\,\d{x}_1 \biggr)\nonumber\\[-8pt]\\[-8pt]
&&\hspace*{40pt}\times\exp \Biggl(-\sum_{i=2}^{n}a_ix_i \Biggr)\,\d
{x}_2\cdots\,\d{x}_{n}\nonumber\\
&=& \frac{f_{n-1}(a_2,\ldots,a_n)-e^{-a_1}f_{n-1}(a_2-a_1,\ldots
,a_n-a_1)}{a_1}.\nonumber
\end{eqnarray}
Alternatively, if we denote $g_0=1$, and for $n\ge1$ and $b_1,\ldots,b_n>0$,
%
\begin{equation}
g_n(b_1,\ldots,b_n)=f_n(b_1,b_1+b_2,\ldots,b_1+\cdots+b_n).
\end{equation}
Then
%
\begin{equation}\hspace*{20pt}
g_n(b_1,\ldots,b_n)=\frac{g_{n-1}(b_1+b_2,b_3,\ldots
,b_n)-e^{-b_1}g_{n-1}(b_2,b_3,\ldots,b_n)}{b_1}.
\end{equation}

From the above, it is also clear (see also \cite{Yor2001}, Theorem~1,
page~31 for the case of a Brownian motion) that, in fact,
%
\begin{eqnarray}
EX^n_t&=&t^nn!f_n(\mu_1t,\ldots,\mu_nt)\nonumber\\[-8pt]\\[-8pt]
&=&t^nn!g_n\bigl(\mu_1t,(\mu_2-\mu
_1)t,\ldots,(\mu_n-\mu_{n-1})t\bigr)\nonumber
\end{eqnarray}
is a linear combination of exponentials. An algorithm for computing the
coefficients of this linear combination is equivalent to the above
simple algorithm which involves only a finite number of additions and
multiplications.

We emphasize that the fact that Theorem~\ref{th:puredeath} holds for
all $n\ge1$, and the algorithm for the computation of moments, also
valid for all $n\ge1$, is special for the case where $Z$ is a nonzero
subordinator. This is true since this is the only case where $\eta
_j(n)$ is finite, strictly positive for all $n\ge1$ and strictly increasing.


\printaddresses

\end{document}